\documentclass[12pt]{amsart}
\usepackage{amsmath, amssymb, latexsym, amsthm,bm}
\usepackage{mathtools}
\usepackage{url}
\usepackage[T1]{fontenc}
\usepackage{mathrsfs}
\usepackage{color}
\usepackage{siunitx}
\usepackage{graphicx}
\usepackage{stmaryrd,amsbsy}
\usepackage{enumitem}
\usepackage{enumerate}
\usepackage{xcolor}
\usepackage{hyperref}
\usepackage{cancel}
\usepackage{xspace}
\usepackage{subcaption}
\usepackage[mathlines]{lineno}
\usepackage{url}
\usepackage{longtable}
\usepackage[obeyFinal,colorinlistoftodos,textsize=tiny]{todonotes}
\usepackage[font=small]{caption}
\usepackage{mathtools,amsmath}
\usepackage{array}
\usepackage{stmaryrd}

\def\genus{{\text{\rm gen}}}
\def\cross#1#2{{\times_{#1,#2}}}

\makeatletter
\def\inmod#1{\allowbreak\mkern5mu({\operator@font mod}\,\,#1)}
\makeatother

\def\edge#1#2{
{#1_{#2}}
}

\def\croge{{\text{\rm crg}}}
\def\cy#1#2{{C_{#1,#2}}}
\def\blueit#1{{\textcolor{blue}{#1}}}
\def\redit#1{{\textcolor{red}{#1}}}

\def\Cross#1#2{{\Scale[2.4]{\cross{#1}{#2}}}}
\def\Tz#1{{\Scale[3.0]{#1}}}

\newcommand*{\Scale}[2][4]{\scalebox{#1}{\ensuremath{#2}}}%

\newtheorem{theorem}{Theorem} 
\newtheorem{theorem*}{Theorem} 
\newtheorem{proposition}[theorem]{Proposition} 

\newtheorem{lemma}[theorem]{Lemma}

\newtheorem{claim}[theorem]{Claim}

\renewcommand{\descriptionlabel}[1]%
  {{\hglue -0.7 cm}\hspace{\labelsep}#1}

\begin{document}



\title[Rotation systems and simple drawings in surfaces]{Rotation systems and simple drawings in surfaces}


\author[Paul]{Rosna Paul}
\address{Institute of Software Technology, Graz University of Technology, Austria}
\email{\tt ropaul@ist.tugraz.at}

\author[Salazar]{Gelasio Salazar}
\address{Instituto de F\'\i sica, Universidad Aut\'onoma de San Luis Potos\'{\i}, SLP 78000, Mexico}
\email{\tt gsalazar@ifisica.uaslp.mx}

\author[Weinberger]{Alexandra Weinberger}
\address{Institute of Software Technology, Graz University of Technology, Austria}
\email{\tt weinberger@ist.tugraz.at}


\makeatletter
\@namedef{subjclassname@2020}{%
  \textup{2020} Mathematics Subject Classification}
\makeatother


\date{\today}

\begin{abstract}
Every simple drawing of a graph in the plane naturally induces a rotation system, but it is easy to exhibit a rotation system that does not arise from a simple drawing in the plane. We extend this to all surfaces: for every fixed surface $\Sigma$, there is a rotation system that does not arise from a simple drawing in $\Sigma$.
\end{abstract}

\maketitle

\section{Introduction}\label{sec:intro}
In a {\em drawing} of a graph, distinct vertices are represented by distinct points in the plane, and each edge $e=uv$ is represented by a Jordan arc whose endpoints are the points that represent $u$ and $v$. It is also required that (D1) no edge contains a vertex other than its endvertices; (D2) every pair of edges intersect each other a finite number of times; (D3) every intersection of edges is either a common endvertex or a crossing (rather than a tangential intersection); and (D4) no three edges cross at a common point. 

A drawing of a graph in an orientable surface determines a {\em rotation} $\pi(v)$ at each vertex $v$: this is the cyclic permutation that records the (clockwise) order in which the edges incident with $v$ leave $v$. The set of the rotations of all the vertices is the {\em rotation system} of the drawing. 

\begin{figure}[ht!]
\centering
\scalebox{0.38}{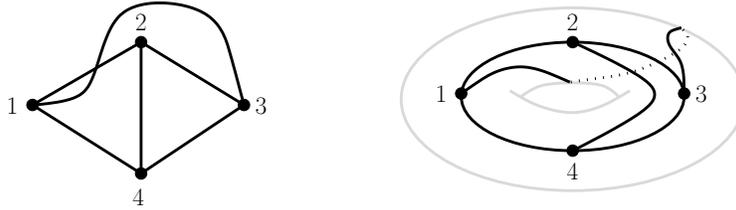}
\caption{Two drawings with the same rotation system $\Pi=\{\pi(1)=(2,3,4),\pi(2)=(1,3,4),\pi(3)=(1,4,2),\pi(4)=(1,2,3)\}$. Note that the drawing on the left-hand side is not simple, since the adjacent edges $12$ and $13$ cross each other. It is easy to verify that this rotation system cannot be realized by a simple drawing in the plane. As illustrated on the right-hand side, $\Pi$ can be realized by a simple drawing in the torus.}
\label{fig:200}
\end{figure}

Rotation systems can be considered independently of their relationship with graph drawings. Let $S$ be a set. A {\em rotation} of an element $s\in S$ is a cyclic permutation $\pi(s)$ of a subset of $S\setminus\{s\}$. A {\em rotation system} on $S$ is a collection $\Pi=\{\pi(s)\}_{s\in S}$ of rotations of all the elements of $S$. If for each $s\in S$ we have that $\pi(s)$ is a cyclic permutation of the entire set $S\setminus\{s\}$, then $\Pi$ is a {\em complete rotation system}. Thus a drawing of a complete graph in an orientable surface determines a complete rotation system.

It is easy to see that every rotation system is the rotation system of a drawing of a graph in the plane; that is, every rotation system can be {\em realized} by a drawing in the plane. However, this does not hold if we slightly restrict the class of drawings under consideration.

One is often interested in drawing a graph so that the total number of edge crossings is as small as possible. It is easy to see that in this context we may focus on {\em simple} drawings, that is, drawings in which (D5) any two edges are either disjoint or have exactly one point in common (which, in view of (D3), is either a common endvertex or a crossing). For instance, the drawing of $K_4$ in the plane on the left-hand side in Figure~\ref{fig:200} is not simple.

Simple drawings thus play a paramount role in topological graph theory, and so it is natural to ask: which rotation systems can be realized by a simple drawing in the plane? Equivalently: which rotation systems are {\em simply realizable} in the plane?
 
As far as we know, this question was first asked by Dan Archdeacon~\cite{Arch1}, for complete rotation systems.
Archdeacon's question was settled by Kyn\v{c}l in~\cite{Kyncl1} (see also~\cite{EuroCG15} and~\cite{Kyncl2}): the decision problem of whether a given complete rotation system is simply realizable in the plane is in {\bf P}. For arbitrary rotation systems little seems to be known, and serious difficulties seem to arise even if one imposes additional restrictions on  the drawings under consideration~\cite{cardinalfelsner}.

As we discuss in more detail in Section~\ref{sec:con}, Archdeacon also considered the simple realizability of rotation systems on compact orientable surfaces. For brevity, throughout this paper we shall refer to a compact orientable surface simply as a {\em surface}. 

When we began to investigate the simple realizability of rotation systems in surfaces, we quickly realized that we needed to start with a question that we had implicitly misjudged as a routine exercise: given a surface $\Sigma$, does there exist a rotation system that cannot be simply realized in $\Sigma$? To our surprise, answering this question was far from trivial.

\begin{theorem}\label{thm:main}

For each surface $\Sigma$, there is a rotation system that is not simply realizable in $\Sigma$. 
\end{theorem}

With the exception of the final discussion in Section~\ref{sec:con}, the rest of the paper is devoted to the proof of Theorem~\ref{thm:main}.

\section{Reducing Theorem~\ref{thm:main} to three lemmas}

To prove Theorem~\ref{thm:main} we define, for each positive integer $n$, a rotation system $\Pi_{n}$ on the set $[n]\cup \{b,r\}$ and prove the following: if $\Sigma$ is a fixed surface, and $n$ is large enough, then $\Pi_{n}$ cannot be realized in $\Sigma$. 

\begin{figure}[ht!]
\centering
\scalebox{0.38}{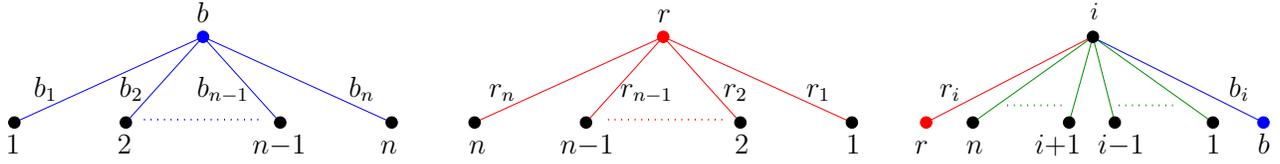}
\caption{Illustration of the rotation system $\Pi_n$.}
\label{fig:300}
\end{figure}

As illustrated in Figure~\ref{fig:300}, the description of $\Pi_n$ is quite simple:
\[  \pi(b) =   \text{\hglue 0.01 cm }(n,\cdots,2,1), \text{\rm \hglue 0.1 cm }    \pi(r) =  \text{\hglue 0.01 cm}(1,2,\ldots,n), \text{\hglue 0.2 cm \rm and \hglue 0.1 cm}   \pi(i) =  \text{\hglue 0.01 cm}(b,1,2,\ldots,i{-}1,i{+}1,\ldots,n,r) \text{\hglue 0.2 cm \rm for } i\in [n] 
\]

\noindent with obvious adjustments for $i=1$ and $i=n$. In a simple drawing that realizes $\Pi_n$, the stars with centers $b,r$, and each $i\in [n]$ are drawn as illustrated in Figure~\ref{fig:300}. 

\def\te#1{{\Scale[2.4]{#1}}}
\def\tf#1{{\Scale[2.4]{#1}}}
\def\tz#1{{\Scale[2.4]{#1}}}
\def\te#1{{\Scale[2.8]{#1}}}
\def\tf#1{{\Scale[2.4]{#1}}}
\def\tz#1{{\Scale[2.0]{#1}}}

To help comprehension, in a drawing that realizes $\Pi_n$ we say that $b$ is the {\em blue} vertex, and the edge $\edge{b}{i}$ that joins $b$ to $i$ is a {\em blue} edge, for $i=1,\ldots,n$. Similarly, vertex $r$ is {\em red}, and the edge $\edge{r}{i}$ that joins $r$ to $i$ is a {\em red} edge, for $i=1,\ldots,n$. Finally, for each distinct $i,j\in[n]$, the edge joining vertices $i$ and $j$ is a {\em green} edge. For each distinct $i,j\in[n]$, the edges $\edge{b}{i},\edge{b}{j},\edge{r}{j},\edge{r}{i}$ form a $4$-cycle, the {\em canonical cycle} $\cy{i}{j}$. 

Note that a drawing that realizes $\Pi_n$ is a drawing of a graph $L_n$ that is isomorphic to $K_{n+2}$ minus one edge, since the only two elements in $[n]\cup \{b,r\}$ that are not adjacent to each other under $\Pi_n$ are $b$ and $r$.

\begin{figure}[ht!]
\centering
\scalebox{0.46}{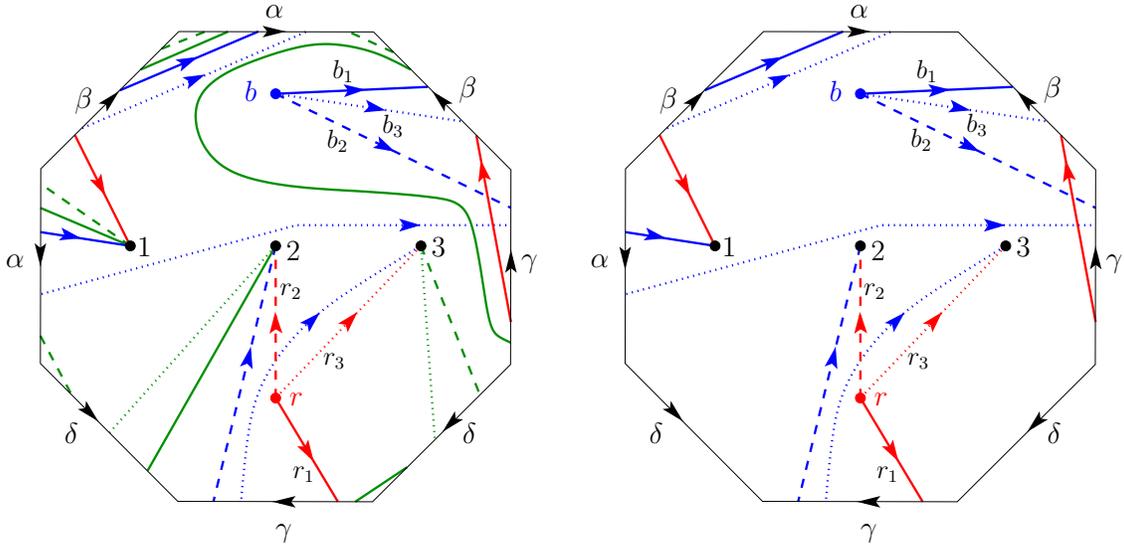}
\caption{On the left-hand side we illustrate a simple drawing of $L_3$ that realizes $\Pi_3$ in the double torus, with three crossings. One crossing involves $\edge{b}{2}$ (dashed) and $\edge{r}{1}$ (solid), another crossing involves $\edge{b}{3}$ (dotted) and $\edge{r}{1}$, and the other crossing involves $\edge{b}{3}$ and $\edge{r}{2}$ (dashed). Thus each canonical cycle has exactly one self-crossing, and so this is a $1$-drawing of $L_3$ in the double torus. If we omit the green edges, we obtain the $1$-drawing of $K_{2,3}$ in the double torus shown on the right-hand side.}
\label{fig:56002}
\end{figure}

It is easy to see that in a simple drawing that realizes $\Pi_n$, each canonical cycle is either crossing-free, or crosses itself exactly once, or crosses itself exactly twice. For $\ell=0,1,2$, we say that an $\ell$-{\em drawing} of $L_n$ is a simple drawing that realizes $\Pi_n$ in some surface, in such a way that every canonical cycle has exactly $\ell$ self-crossings. In Figure~\ref{fig:56002} we illustrate a $1$-drawing of $L_3$ in the double torus. We use the standard polygonal representation of the double torus as an octagon whose sides are identified in pairs, according to the indicated labels and orientations.

\def\tf#1{{\Scale[2.0]{#1}}}
\def\tz#1{{\Scale[1.8]{#1}}}
\def\tf#1{{\Scale[2.4]{#1}}}
\def\tz#1{{\Scale[2.0]{#1}}}

If we discard the green edges from a drawing $D_n$ that realizes $\Pi_n$, the result is a drawing $E_n$ of $K_{2,n}$ in which the rotations of the vertices are $\pi(b) = (n,\cdots,2,1), \pi(r) = (1,2,\ldots,n)$, and $\pi(i) = (b,r)$, for $i\in [n]$. Note that the notion of a canonical cycle carries over naturally to $E_n$. If $D_n$ is an $\ell$-drawing for some $\ell\in\{0,1,2\}$, then clearly each canonical cycle also has $\ell$ self-crossings in $E_n$, and we also say that $E_n$ is an {\em $\ell$-drawing} of $K_{2,n}$. As it happens, most of the proof of Theorem~\ref{thm:main} involves $1$-drawings and $2$-drawings of $K_{2,n}$.

The following proposition is crucial, as it will allow us to focus for the rest of the paper on drawings of $L_n$ that are $\ell$-drawings for some $\ell\in\{0,1,2\}$.

\begin{proposition}\label{pro:el0}
Let $\Sigma$ be a surface. Suppose that for each positive integer $m$, there is a simple drawing of $L_m$ that realizes $\Pi_m$ in $\Sigma$. Then for each positive integer $n$, there is an $\ell$-drawing of $L_n$ in $\Sigma$, for some $\ell\in\{0,1,2\}$.
\end{proposition}

\begin{proof}
Let $n$ be a positive integer, and let $m$ be the Ramsey number $R_2(n,n,n)$: if we colour each edge of the complete graph $K_m$ with one of three available colours, then there is a complete subgraph of size $n$ all of whose edges are of the same colour.

By hypothesis, there is a simple drawing $D_m$ that realizes $\Pi_m$ in $\Sigma$. We construct an auxiliary complete graph $G$ with vertex set $[m]$, and colour an edge $ij$ of $G$ with colour $\ell\in\{0,1,2\}$ if the canonical cycle $\cy{i}{j}$ has exactly $\ell$ self-crossings in $D_m$. Since each canonical cycle has either $0,1$, or $2$ self-crossings, this edge colouring of $G$ is well-defined.

By Ramsey's theorem there is a complete subgraph $H$ of $G$ of size $n$, all of whose edges are of the same colour $\ell\in\{0,1,2\}$. Let $i_1,i_2,\ldots,i_n$ be the vertices of $H$, labelled so that $i_1 < i_2 < \cdots < i_n$. Delete from $D_m$ all the vertices in $[m]\setminus\{i_1,i_2,\ldots,i_n\}$ and their incident edges, and for each $j=1,\ldots,n$ relabel the vertex $i_j$ with $j$. Clearly, as a result we obtain an $\ell$-drawing of $L_n$ in $\Sigma$.
\end{proof}

In view of Proposition~\ref{pro:el0}, to prove Theorem~\ref{thm:main} it suffices to show that, for each $\ell\in\{0,1,2\}$, no fixed surface can host $\ell$-drawings of $L_n$ for arbitrarily large values of $n$. This is achieved with the next three lemmas, which take care separately of the three possible values of $\ell$. Actually, in Lemmas~\ref{lem:1cr} and~\ref{lem:2cr} we establish results that are strictly stronger than we need to prove Theorem~\ref{thm:main}, as they involve $1$- and $2$-drawings of $K_{2,n}$, a very sparse subgraph of $L_n$.

\begin{lemma}\label{lem:0cr}
For each surface $\Sigma$, there is an $n$ such that there is no $0$-drawing of $L_n$ in $\Sigma$.
\end{lemma}

\begin{lemma}\label{lem:1cr}
For each surface $\Sigma$, there is an $n$ such that there is no $1$-drawing of $K_{2,n}$ in $\Sigma$.
\end{lemma}

\begin{lemma}\label{lem:2cr}
For each surface $\Sigma$, there is an $n$ such that there is no $2$-drawing of $K_{2,n}$ in $\Sigma$.
\end{lemma}

The rest of the paper is devoted to the proofs of these lemmas. Assuming the lemmas, Theorem~\ref{thm:main} immediately follows. We include the proof for completeness.

\begin{proof}[Proof of Theorem~\ref{thm:main}, assuming Lemmas~\ref{lem:0cr},~\ref{lem:1cr}, and~\ref{lem:2cr}] Let $\Sigma$ be a surface. Note that Lemma~\ref{lem:1cr} (respectively, Lemma~\ref{lem:2cr}) implies that there is an $n$ such that there is no $1$-drawing (respectively, $2$-drawing) of $L_n$ in $\Sigma$. Combining this with Lemma~\ref{lem:0cr}, we conclude that there is an $N$ such that there is no $\ell$-drawing of $L_N$ in $\Sigma$, for any $\ell\in\{0,1,2\}$. 

The contrapositive of Proposition~\ref{pro:el0} then implies that there is an integer $m$ such that there is no simple drawing that realizes $\Pi_m$ in $\Sigma$. In other words, $\Pi_m$ is not simply realizable in $\Sigma$.
\end{proof}

\section{Proof of Lemma~\ref{lem:0cr}}\label{sec:proof0}

To prove Lemma~\ref{lem:0cr}, we establish the next two propositions, where we use the following terminology. A {\em main cycle} of $L_n$ is a canonical cycle of the form $\cy{i}{i+1}$, for some $i=1,\ldots,n$ (indices are read modulo $n$). In a $0$-drawing of $L_n$ in some surface, a main cycle is {\em good} if it bounds a disk. Otherwise, it is {\em bad}. We note that the notions of a main cycle, of a good cycle, and of a bad cycle, carry over naturally to $0$-drawings of $K_{2,n}$. 

\begin{proposition}\label{pro:01}
Let $n\ge 3$ be an integer. If $D_n$ is a $0$-drawing of $L_n$ in a surface, then every main cycle is bad in $D_n$.
\end{proposition}

In the next statement, and in the rest of the paper, we use $\genus(\Sigma)$ to denote the genus of a surface $\Sigma$. 

\begin{proposition}\label{pro:02}
Let $n\ge 3$ be an integer. If $E_n$ is a $0$-drawing of $K_{2,n}$ in a surface $\Sigma$, then at least $n - 2\,\genus(\Sigma)$ main cycles are good in $E_n$.
\end{proposition}

Deferring the proofs of these statements for a moment, we note that Lemma~\ref{lem:0cr} follows easily from them.

\begin{proof}[Proof of Lemma~\ref{lem:0cr}]
If $\Sigma$ has genus $0$ we let $n=3$. Otherwise, let $n=2\,\genus(\Sigma)+1$. By way of contradiction, suppose there is a $0$-drawing $D_n$ of $L_n$ in $\Sigma$. Note that $D_n$ contains a $0$-drawing $E_n$ of $K_{2,n}$ in $\Sigma$. Thus by Proposition~\ref{pro:02} it follows that there is at least one good main cycle in $E_n$, and so there is at least one good main cycle in $D_n$. But this contradicts Proposition~\ref{pro:01}.
\end{proof}

We conclude the section with the proofs of Propositions~\ref{pro:01} and~\ref{pro:02}.

\begin{proof}[Proof of Proposition~\ref{pro:01}]
Let $D_n$ be a $0$-drawing of $L_n$ in a surface. By way of contradiction, suppose that the main cycle $\cy{i}{i+1}$ is good in $D_n$, for some $i\in[n]$. Without loss of generality we may then assume that as we traverse in $D_n$ the cycle $\cy{i}{i+1}$ visiting the vertices $b,i,r,i+1$ in this order, we find that on the left-hand side of this traversal we bound a disk $\Delta$. See Figure~\ref{fig:400}. 

\def\tf#1{{\Scale[2.4]{#1}}}
\def\tz#1{{\Scale[2.4]{#1}}}
\begin{figure}[ht!]
\centering
\scalebox{0.4}{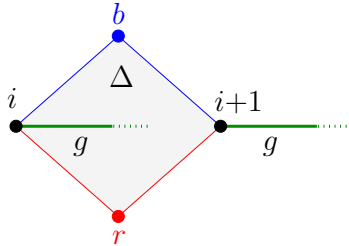}
\caption{Illustration of the proof of Proposition~\ref{pro:01}.}
\label{fig:400}
\end{figure}
\def\tf#1{{\Scale[2.4]{#1}}}
\def\tz#1{{\Scale[2.0]{#1}}}

The rotation at $i$ in $\Pi_n$ contains the cyclic subpermutation $(b,i+1,r)$, and this means that at least part of the green edge $g$ that joins $i$ and $i+1$ must lie inside $\Delta$, as illustrated in Figure~\ref{fig:400}. On the other hand, the rotation at $i+1$ in $\Pi_n$ contains the cyclic subpermutation $(b,i,r)$, and so part of $g$ must lie outside $\Delta$, as we also illustrate in Figure~\ref{fig:400}. By the Jordan curve theorem, it follows that the edge $g$ must cross the cycle $\cy{i}{i+1}$. But this contradicts that $D_n$ is a simple drawing, since $g$ joins $i$ and $i+1$ and each edge of $\cy{i}{i+1}$ is incident with either $i$ or $i+1$.
\end{proof}


\begin{proof}[Proof of Proposition~\ref{pro:02}]
We prove the proposition by induction on the genus of $\Sigma$. In the base case $\genus(\Sigma)=0$, and so $\Sigma$ is the sphere. In this case the proposition claims that all $n$ main cycles in $D_n$ are good. This is trivially true, since the host surface is the sphere. 

The inductive hypothesis is that the proposition holds for surfaces of genus $k-1$, for some integer $k\ge 1$. For the inductive step, let $E_n$ be a $0$-drawing of $K_{2,n}$ in a surface $\Sigma$ of genus $\genus(\Sigma)= k$. Now $E_n$ has $n$ main cycles, and so we assume that $n > 2\,\genus(\Sigma)$, as otherwise there is nothing to prove.  
Note that since $E_n$ is a $0$-drawing of $K_{2,n}$, we have in particular that no edge of $K_{2,n}$ is crossed in $E_n$.

If all main cycles in $E_n$ are good we are done, and so we suppose that there is a main cycle $\cy{i}{i+1}$ in $E_n$ that is bad. Thus $\cy{i}{i+1}$ is a noncontractible Jordan curve in $\Sigma$. As illustrated in Figure~\ref{fig:480}, we take two simple closed curves $\gamma_1,\gamma_2$ homotopic to $\cy{i}{i+1}$, drawn very close to $\cy{i}{i+1}$ so that $E_n$ does not intersect the closed annulus bounded by $\gamma_1$ and $\gamma_2$. Note that such a closed annulus (disjoint from $E_n$) must exist, as otherwise some edge of $\cy{i}{i+1}$ would be crossed in $E_n$.

\begin{figure}[ht!]
\centering
\scalebox{0.55}{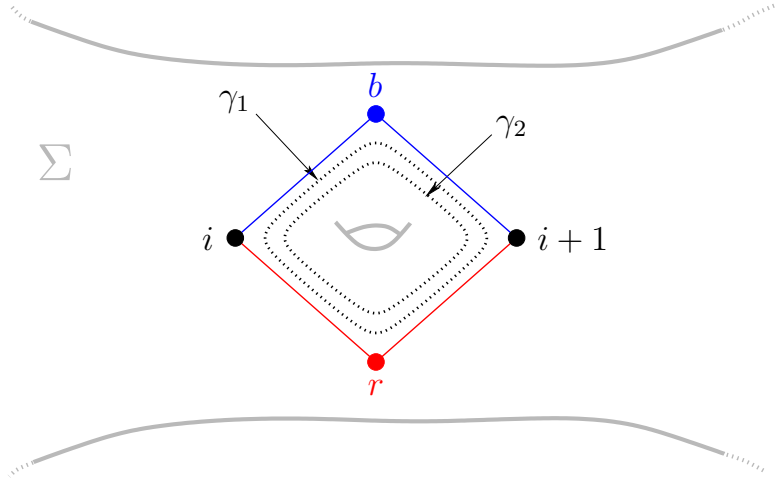}
\caption{Illustration of the proof of Proposition~\ref{pro:02}.}
\label{fig:480}
\end{figure}

We now remove the open annulus bounded by $\gamma_1$ and $\gamma_2$. This turns $\Sigma$ into a compact surface with two boundary components, namely $\gamma_1$ and $\gamma_2$. We glue the boundary of a disk $\Delta_1$ to $\gamma_1$, and we glue the boundary of a disk $\Delta_2$ to $\gamma_2$. As a result, we obtain a (compact, without boundary) surface $\Sigma'$ of genus $\genus(\Sigma)-1 = k- 1$.

The drawing $E_n$ naturally induces a drawing $E_n'$ in $\Sigma'$ with the same properties as $E_n$. That is, $E_n'$ is a $0$-drawing of $K_{2,n}$ in $\Sigma'$. Thus by the induction hypothesis $E_n'$ has at least $n-2\,\genus(\Sigma')=n-2(\,\genus(\Sigma)-1)=n-2\,\genus(\Sigma)+2$ main cycles that are good. To finish the proof, we show that at most two of these main cycles are bad in $E_n$. 

Note that the main cycle $\cy{i}{i+1}$ is bad in $E_n$, but it became good in $E_n'$, as it now bounds a disk that contains $\Delta_1$. The key point is that there is at most one other main cycle that is bad in $E_n$ but became good in $E_n'$: this would be a main cycle that in $\Sigma$ was homotopic to $\cy{i}{i+1}$ (and hence it was bad in $E_n$), but became good in $E_n'$ because in $\Sigma'$ it bounds a disk that contains~$\Delta_2$.

From this discussion, we conclude that the number of main cycles that are good in $E_n$ is at least the number of main cycles that are good in $E_n'$ minus two. Since this last number is at least $n-2\,\genus(\Sigma)+2$, it follows that the number of main cycles that are good in $E_n$ is at least $n-2\,\genus(\Sigma)$.
\end{proof}

\section{Narrowing down the class of $1$-drawings we need to consider}\label{sec:narrow1}

At a high level, the strategy of the proof of Lemma~\ref{lem:1cr} has two main steps. The first step, which is the goal of this section, is the following. Let $\Sigma$ be a surface. Seeking a contradiction, we show that if for every $n$ there is a $1$-drawing of $K_{2,n}$ in $\Sigma$, then there is one such drawing that satisfies certain highly restrictive properties (to be specified shortly). The second step, achieved in the next section, will consist of showing that no fixed surface can host arbitrarily large $1$-drawings that satisfy these properties, thus providing the desired contradiction.

Throughout this section and Section~\ref{sec:1cr}, whenever we refer to a $1$-drawing we mean a $1$-drawing of $K_{2,n}$ in a surface $\Sigma$, for some positive integer $n$. We use several times expressions such as ``if there are arbitrarily large $1$-drawings in $\Sigma$''. This is to be interpreted as ``if there exist $1$-drawings of $K_{2,n}$ in $\Sigma$ for arbitrarily large values of $n$''.

\begin{figure}[ht!]
\centering
\scalebox{0.45}{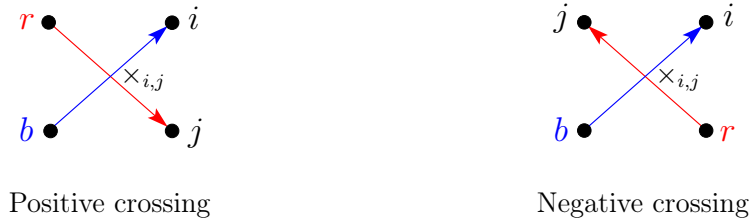}
\caption{Every crossing in a simple drawing of $K_{2,n}$ is positive or negative, according to this convention.}
\label{fig:800}
\end{figure}

\begin{figure}[ht!]
\centering
\scalebox{0.44}{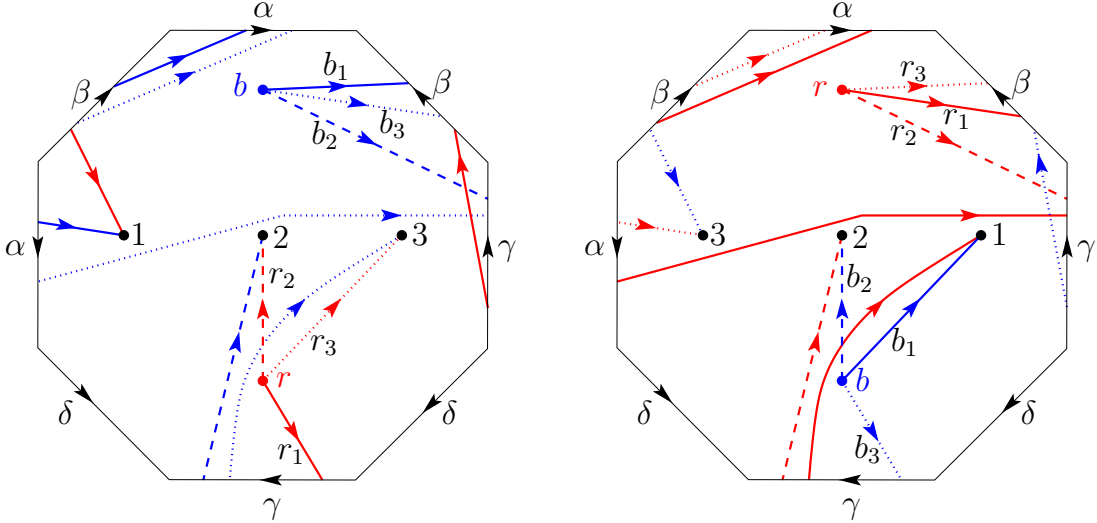}
\caption{On the left-hand side we illustrate a $1$-drawing of $K_{2,3}$ in the double torus. Each crossing in this drawing is negative, and so this $1$-drawing is negative. The three crossings $\cross{2}{1}, \cross{3}{1},\cross{3}{2}$ in this drawing are backward crossings, and so this is a backward $1$-drawing. If we perform the relabellings $b\mapsto r, r\mapsto b, 1\mapsto 3$, and $3\mapsto 1$, we obtain the $1$-drawing on the right-hand side, which is also backward, but positive. 
}
\label{fig:5400}
\end{figure}

Let $D_n$ be a $1$-drawing of $K_{2,n}$. Since $D_n$ is simple, each crossing in $D_n$ involves a blue edge~$\edge{b}{i}$ and a red edge $\edge{r}{j}$, for some distinct $i,j\in[n]$. We use $\cross{i}{j}$ to denote this crossing. We assign a sign to $\cross{i}{j}$ as follows. Orient each blue edge $\edge{b}{i}$ from $b$ to $i$, and orient each red edge $\edge{r}{j}$ from $r$ to $j$. As illustrated in Figure~\ref{fig:800}, a crossing is {\em positive} if the red edge crosses the blue edge from the left-hand side of the blue edge. Otherwise, the crossing is {\em negative}.  
If all the crossings in $D_n$ are positive (respectively, negative), then we say that $D_n$ itself is {\em positive} (respectively, {\em negative}). 
On the left-hand side of Figure~\ref{fig:5400} we illustrate a negative $1$-drawing of $K_{2,3}$, and on the right-hand side we show a positive $1$-drawing of $K_{2,3}$.

%
%

\def\tf#1{{\Scale[2.4]{#1}}}
\def\tz#1{{\Scale[2.4]{#1}}}
\def\tf#1{{\Scale[2.4]{#1}}}
\def\tz#1{{\Scale[2.0]{#1}}}

We start to narrow down the $1$-drawings of interest with the following statement.

\begin{proposition}\label{pro:1cr1}
Let $\Sigma$ be a surface. If there exist arbitrarily large $1$-drawings in $\Sigma$, then there exist arbitrarily large $1$-drawings in $\Sigma$ that are positive.
\end{proposition}

\begin{proof}
Let $n$ be any positive integer. Suppose that {there exist arbitrarily large $1$-drawings in $\Sigma$}. Similarly as in the proof of Proposition~\ref{pro:el0}, an application of Ramsey's theorem yields that there exists a $1$-drawing $D_n$ of $K_{2,n}$ in $\Sigma$ that is either positive or negative.

If $D_n$ is positive, then we are done. If it is negative, we proceed as illustrated in Figure~\ref{fig:5400}: we relabel the vertices $b \mapsto r, r\mapsto b$, and $i\mapsto n-i+1$ for $i=1,\ldots,n$, and as a result we obtain a $1$-drawing of $K_{2,n}$ that is positive. 
\end{proof}

In order to continue the process of narrowing down the $1$-drawings of interest, we introduce the notions of forward and backward crossings. 
Let $\cross{i}{j}$ be a crossing in a $1$-drawing $D_n$ of $K_{2,n}$ in some surface. Recall that this is a crossing between the blue edge $\edge{b}{i}$ and the red edge~$\edge{r}{j}$, where~$i\neq j$. If $i<j$, then $\cross{i}{j}$ is a {\em forward} crossing, and if $i>j$ then it is a {\em backward}~crossing. If all the crossings are forward (respectively, backward) then $D_n$ itself is a {\em forward} (respectively, {\em backward}) drawing. 
For instance, both $1$-drawings of $K_{2,3}$ shown in Figure~\ref{fig:5400} are backward drawings.

\begin{proposition}\label{pro:1cr2}
Let $\Sigma$ be a surface. If there exist arbitrarily large $1$-drawings in $\Sigma$ that are positive, then there exist arbitrarily large $1$-drawings in $\Sigma$ that are positive and forward.
\end{proposition}

\begin{figure}[ht!]
\centering
\scalebox{0.4}{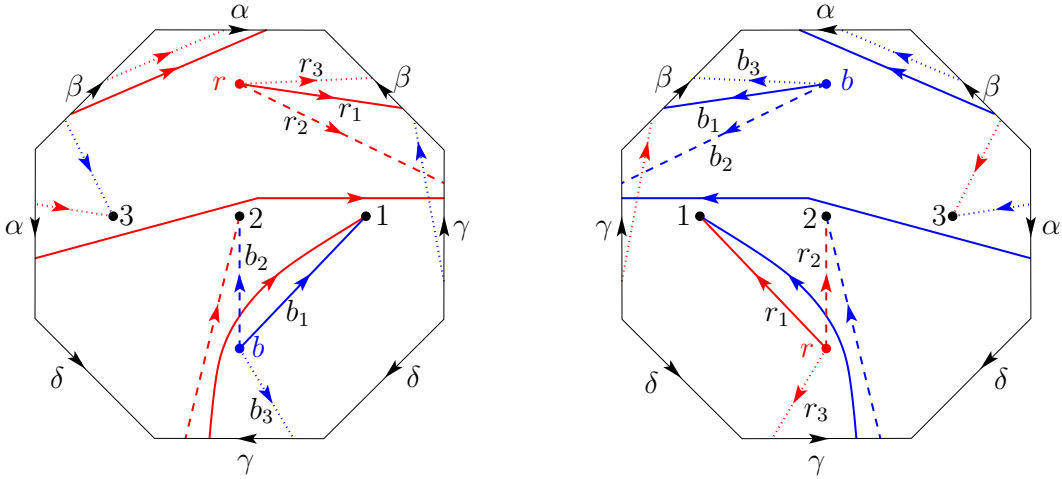}
\caption{On the left-hand side we have a $1$-drawing of $K_{2,3}$, in the double torus, that is positive and backward. If we perform a reflection on the polygon (that is, an orientation-reversing self-homeomorphism of the double torus) and exchange the labels of $b$ and $r$, we obtain the drawing on the right-hand side. This is a $1$-drawing of $K_{2,3}$ that is positive and forward.}
\label{fig:5800}
\end{figure}

\begin{proof}
Let $n$ be a positive integer. Suppose that {there exist arbitrarily large positive $1$-drawings in $\Sigma$}. Similarly as in the proof of Proposition~\ref{pro:el0}, an application of Ramsey's theorem yields that there exists a $1$-drawing $D_n$ of $K_{2,n}$ in $\Sigma$ that is positive and either forward or backward.

If $D_n$ is forward then we are done. If it is backward, we transform it as illustrated in Figure~\ref{fig:5800}. First, we perform an orientation-reversing self-homeomorphism on $\Sigma$, and let $D_n'$ denote the drawing that is the image of $D_n$ under this mapping. Finally, we modify $D_n'$ by exchanging the labels of $b$ and $r$. The result is a $1$-drawing of $K_{2,n}$ that is positive and forward.
\end{proof}

\def\tf#1{{\Scale[2.4]{#1}}}
\def\tz#1{{\Scale[2.4]{#1}}}
\def\tf#1{{\Scale[2.4]{#1}}}
\def\tz#1{{\Scale[2.0]{#1}}}

To complete this process of narrowing down the $1$-drawings that we need to consider, we introduce one final notion. Let $D_n$ be a $1$-drawing of $K_{2,n}$. As we traverse a blue edge $\edge{b}{i}$ following its direction (we recall, from $b$ to $i$), we encounter its crossings with red edges in some order $\cross{i}{j_1},\cross{i}{j_2},\ldots,\cross{i}{j_s}$. We say that $j_1,j_2,\ldots,j_s$ is the {\em crossing sequence} of $\edge{b}{i}$.

If the crossing sequence of $\edge{b}{i}$ is increasing (respectively, decreasing), then we say that the edge $\edge{b}{i}$ is {\em increasing} (respectively, {\em decreasing}). If every blue edge is increasing (respectively, decreasing), then $D_n$ is {\em blue increasing} (respectively, {\em blue decreasing}).

Similarly, as we traverse a red edge $\edge{r}{j}$ following its direction (from $r$ to $j$) we encounter its crossings with blue edges in some order $\cross{i_1}{j},\cross{i_2}{j},\ldots,\cross{i_t}{j}$. Then $i_1,i_2,\ldots,i_t$ is the {\em crossing sequence} of $\edge{r}{i}$. If this sequence is increasing (respectively, decreasing), then $\edge{r}{j}$ is {\em increasing} (respectively, {\em decreasing}). If every red edge is increasing (respectively, decreasing), then $D_n$ is {\em red increasing} (respectively, {\em red decreasing}).

If $D_n$ is blue increasing or blue decreasing, then it is {\em blue monotone}. If $D_n$ is red increasing or red decreasing, then it is {\em red monotone}. Finally, if $D_n$ is both blue monotone and red monotone, then it is {\em monotone}.

\begin{proposition}\label{pro:1cr3}
Let $\Sigma$ be a surface. If there exist arbitrarily large $1$-drawings in $\Sigma$ that are positive and forward, then there exist arbitrarily large $1$-drawings in $\Sigma$ that are positive, forward, and monotone.
\end{proposition}

\begin{proof}
We prove that if ($*$) {\sl there exist arbitrarily large $1$-drawings in $\Sigma$ that are positive and forward}, then ($**$) {\sl there exist arbitrarily large $1$-drawings in $\Sigma$ that are positive, forward, and blue monotone.} A totally analogous sequence of arguments shows that ($**$) implies that there exist arbitrarily large $1$-drawings in $\Sigma$ that are positive, forward, blue monotone, and red monotone (that is, positive, forward, and monotone, as claimed in the proposition).

Let $n$ be any positive integer. Assuming ($*$), our aim is to show that there is a $1$-drawing of $K_{2,n}$ in $\Sigma$ that is positive, forward, and blue monotone.

Let $m$ be the Ramsey number $R_3(n,n)$: if all the $3$-edges in a complete $3$-uniform hypergraph of size $m$ are coloured with one of two available colours, then there is a complete subhypergraph of size $n$ all of whose $3$-edges are of the same colour.

By assumption, there is a $1$-drawing $D_m$ of $K_{2,m}$ in $\Sigma$ that is positive and forward. Construct an auxiliary complete $3$-uniform hypergraph $G$, whose vertex set is $[m]$. Let $e=\{p,q,s\}$ be a $3$-edge of $G$, where $1 {\le} p{<}q{<}s \le m$. Note that since $D_m$ is forward, the blue edge $b_p$ is crossed by the red edges $r_{q}$ and $r_s$. We colour $e$ black if as we traverse the blue edge $b_p$ in $D_m$ from $b$ to~$p$, we encounter its crossing with $r_q$ before its crossing with $r_s$, and we colour $e$ green otherwise.

Clearly, every $3$-edge of $G$ is either black or green, and so by Ramsey's theorem there exist integers $1 \le i_1 < i_2 < \cdots < i_n \le m$ such that in the complete subhypergraph $H$ of $G$ on $\{i_1,i_2,\ldots,i_n\}$, all the $3$-edges are of the same colour. Remove from $D_m$ all the vertices in $[m]\setminus\{i_1,i_2,\ldots,i_n\}$ and their incident edges, and relabel the remaining $n$ vertices with the rule $i_j\mapsto j$ for $j=1,\ldots,n$. Let $D_n$ be the resulting $1$-drawing of $K_{2,n}$.

If all the $3$-edges of $H$ are black, then $D_n$ is blue increasing, and if all the $3$-edges of $H$ are green then $D_n$ is blue decreasing. That is, $D_n$ is blue monotone. We finish the proof by noting that since $D_m$ is positive and forward, then clearly $D_n$ is also positive and forward.
\end{proof}

Combining Propositions~\ref{pro:1cr1},~\ref{pro:1cr2}, and~\ref{pro:1cr3}, we obtain the following statement. This proposition will allow us, in the proof of Lemma~\ref{lem:1cr}, to focus our attention on $1$-drawings that are positive, forward, and monotone.

\begin{proposition}\label{pro:1cr4}
Let $\Sigma$ be a surface. If there exist arbitrarily large $1$-drawings in $\Sigma$, then there exist arbitrarily large $1$-drawings in $\Sigma$ that are positive, forward, and monotone.
\end{proposition}

\section{Proof of Lemma~\ref{lem:1cr}}\label{sec:1cr}

Lemma~\ref{lem:1cr} claims that no fixed surface can host $1$-drawings of $K_{2,m}$ for arbitrarily large values of $m$. In view of Proposition~\ref{pro:1cr4}, in order to prove the lemma it suffices to show that {\em no fixed surface can host $1$-drawings of $K_{2,n}$ that are positive, forward, and monotone, for arbitrarily large values of $n$.}

As we shall see shortly, Lemma~\ref{lem:1cr} in this form follows easily from the next two claims. 

\begin{claim}\label{cla:1cla1}
Let $D_n$ be a $1$-drawing of $K_{2,n}$ in a surface $\Sigma$, where $n$ is even. Suppose that $D_n$ is positive and forward. If $D_n$ is blue increasing or red increasing, then $\genus(\Sigma)\ge (n-2)/2$.
\end{claim}

\begin{claim}\label{cla:1cla2}
Let $D_n$ be a $1$-drawing of $K_{2,n}$ in a surface $\Sigma$. Suppose that $D_n$ is positive and forward. If $D_n$ is blue decreasing and red decreasing, then $\genus(\Sigma)\ge n-2$.
\end{claim}

\begin{proof}[Proof of Lemma~\ref{lem:1cr}, assuming Claims~\ref{cla:1cla1} and~\ref{cla:1cla2}]
Let $\Sigma$ be a surface. By way of contradiction, suppose that for every $m$ there is a $1$-drawing of $K_{2,m}$ in $\Sigma$. By Proposition~\ref{pro:1cr4}, it follows that ($\dag$) {\sl for every $n$ there is a $1$-drawing of $K_{2,n}$ in $\Sigma$ that is forward, positive, and monotone.}

Let $n > 2+2 \,\genus(\Sigma)$ be an even integer. By ($\dag$), there exists a $1$-drawing $D_n$ of $K_{2,n}$ in $\Sigma$ that is forward, positive, and monotone. Thus $D_n$ is forward, positive, and either (i) blue increasing or red increasing; or (ii) both blue and red decreasing. If (i) holds, we have a contradiction to Claim~\ref{cla:1cla1}, and if (ii) holds, we have a contradiction to Claim~\ref{cla:1cla2}.
\end{proof}

To prove the claims, we make use of the common device of regarding a drawing as an embedding, by turning each crossing into a degree $4$ vertex, a {\em crossing vertex} (coloured white in the figures). Under this perspective, we can use the topological graph theory machinery to investigate drawings. 

Note that under this perspective, what used to be a blue edge $\edge{b}{i}$ in a $1$-drawing $D_n$ of $K_{2,n}$ becomes a blue {path}, as the crossings of $\edge{b}{i}$ are now regarded as vertices, which are the internal vertices of the path $\edge{b}{i}$. Similarly, what used to be a red edge $\edge{r}{i}$ becomes a red path. 

We make extensive use of the following fact (see for instance~\cite[Section 4.1]{moharthomassen}). If a graph with $\nu$ vertices and $\varepsilon$ edges is embedded in $\Sigma$, and the embedding has $\phi$ facial walks, then
\begin{equation}\label{eq:genuseq}
\genus(\Sigma) \ge \frac{1}{2}\biggl(2 - \nu + \varepsilon - \phi\biggr).
\end{equation}

\begin{proof}[Proof of Claim~\ref{cla:1cla1}]
Suppose first that $D_n$ is blue increasing, and recall that by assumption $D_n$ is positive and forward. Let $J_n$ be the restriction of $D_n$ to the blue path $\edge{b}{1}$ and the red paths $\edge{r}{2},\ldots,\edge{r}{n}$ (we suppress the crossing vertices of these red paths with other blue edges). As illustrated in Figure~\ref{fig:4700}, the vertices of $J_n$ are $r,b,1,\ldots,n$, and $\cross{1}{2},\ldots,\cross{1}{n}$, which are the crossings of $\edge{b}{1}$ with the red paths $\edge{r}{2},\ldots,\edge{r}{n}$, respectively.

\def\tf#1{{\Scale[2.55]{#1}}}
\def\tz#1{{\Scale[2.55]{#1}}}
\begin{figure}[ht!]
\centering
\hglue 0cm \scalebox{0.4}{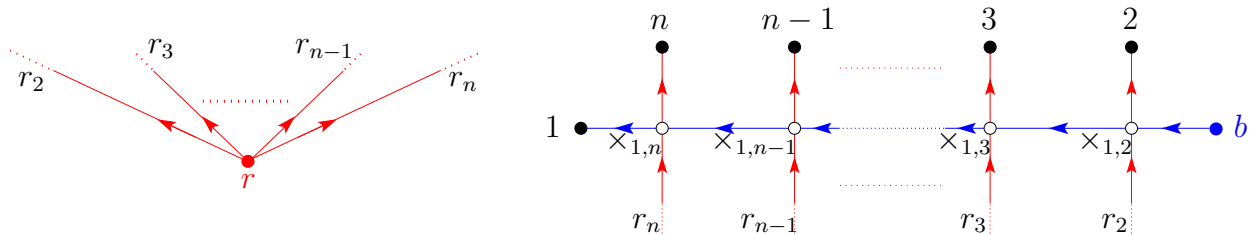}
\caption{Illustration of the proof of Claim~\ref{cla:1cla1}.}
\label{fig:4700}
\end{figure}
\def\tf#1{{\Scale[2.4]{#1}}}
\def\tz#1{{\Scale[2.0]{#1}}}

Thus $J_n$ has $2n+1$ vertices and $3n-2$ edges, and it is easy to verify that it has only one facial walk. (Here is where we use that $n$ is even; if $n$ were odd, there would be two facial walks). Since $\Sigma$ is the host surface of $J_n$, it follows from~\eqref{eq:genuseq} that $\genus(\Sigma)\ge (n-2)/2$.

If $D_n$ is red increasing, we consider instead the restriction $J_n$ of $D_n$ to the red path $\edge{r}{n}$ and the blue paths $\edge{b}{1},\ldots,\edge{b}{n-1}$. As in the previous case, it is easy to verify that $J_n$ has $2n+1$ vertices, $3n-2$ edges, and one facial walk, thus implying that $\genus(\Sigma)\ge (n-2)/2$.
\end{proof}

\begin{proof}[Proof of Claim~\ref{cla:1cla2}]
As in the proof of Claim~\ref{cla:1cla1}, the strategy is to find a restriction $J_n$ of $D_n$ that has exactly one face, and to apply~\ref{eq:genuseq} to obtain the required bound on $\genus(\Sigma)$. As we shall see, the restriction $J_n$ is as illustrated in Figure~\ref{fig:3700}.

\begin{figure}[ht!]
\centering
\hglue -1cm \scalebox{0.38}{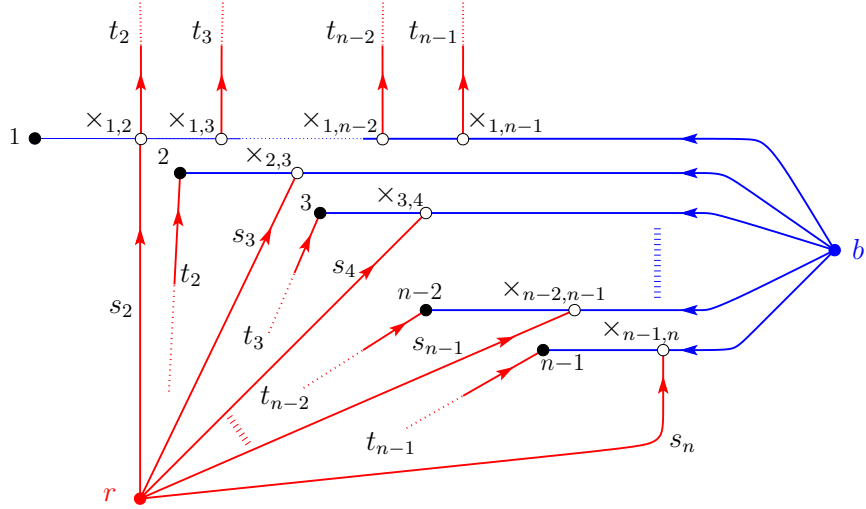}
\caption{Illustration of the proof of Claim~\ref{cla:1cla2}.}
\label{fig:3700}
\end{figure}

Note that since $D_n$ is positive, forward, and red decreasing, for $i=2,\ldots,n$ the red path $\edge{r}{i}$ crosses the blue paths $\edge{b}{i-1},\edge{b}{i-2},\ldots,\edge{b}{1}$ in this order, from their left-hand sides. In particular, as we illustrate in Figure~\ref{fig:3700}, the first crossing of $\edge{r}{i}$ as we traverse it from $r$ to $i$ is $\cross{i-1}{i}$, and the last crossing is $\cross{1}{i}$. 

We also note that since $D_n$ is blue decreasing, as we traverse $\edge{b}{1}$ from $b$ to $1$ we encounter the crossings $\cross{1}{n-1},\cross{1}{n-2},\ldots,\cross{1}{3},\cross{1}{2}$ in this order. See Figure~\ref{fig:3700}.

In the restriction $J_n$ we include the blue paths $\edge{b}{1},\ldots,\edge{b}{n-1}$, and for $i=2,\ldots,n-1$, we include two portions of the red path $\edge{r}{i}$: the part $s_i$ from $r$ to its first crossing $\cross{i-1}{i}$, and the part $t_i$ from its last crossing $\cross{1}{i}$ to vertex $i$. Finally, include from $\edge{r}{n}$ its part $s_n$ from $r$ to its first crossing $\cross{n-1}{n}$. See Figure~\ref{fig:3700}. Thus $J_n$ has $3n{-}3$ vertices and $5n{-}8$ edges, and it is easy to check that it has only one facial walk. Thus it follows from~\eqref{eq:genuseq} that $\genus(\Sigma)\ge n-2$.
\end{proof}
\newpage

\section{Narrowing down the class of $2$-drawings we need to consider}\label{sec:narrow2}

As we did for the proof of Lemma~\ref{lem:1cr}, we pave the way towards the proof of Lemma~\ref{lem:2cr} by showing that it suffices to prove the lemma for a highly restricted class of $2$-drawings of $K_{2,n}$.

Throughout this section and Section~\ref{sec:2cr}, whenever we refer to a $2$-drawing we mean a $2$-drawing of $K_{2,n}$ in a surface $\Sigma$, for some positive integer $n$. We use several times expressions such as ``if there are arbitrarily large $2$-drawings in $\Sigma$''. This is to be interpreted as ``if there exist $2$-drawings of $K_{2,n}$ in $\Sigma$ for arbitrarily large values of $n$''.

We have an additional important remark. In the previous section, for the proof of Lemma~\ref{lem:1cr}, it was convenient to turn crossings into degree $4$ vertices: by regarding drawings as embeddings, we were able to use~\eqref{eq:genuseq}. This device will also be used in the next section for the proof of Lemma~\ref{lem:2cr}. However, in the current section (similarly to Section~\ref{sec:narrow1}) we are interested in properties of {$2$-drawings} of $K_{2,n}$, and so we do {\em not} regard crossings as degree $4$ vertices. Thus for $i=1,\ldots,n$, throughout this section $\edge{b}{i}$ and $\edge{r}{i}$ maintain their identities as edges. 

We start by noting that in a $2$-drawing $D_n$ of $K_{2,n}$, for each $i\in[n]$ the red edge $\edge{r}{i}$ crosses all the blue edges except for $\edge{b}{i}$. Thus the crossing sequence of $\edge{r}{i}$ is a permutation of $[n]\setminus\{i\}$. We say that $\edge{r}{i}$ is {\em ahead first} if in this sequence $i+1,i+2,\ldots,n$ appear first in some order, followed by $1,2,\ldots,i-1$ in some order. If every red edge in $D_n$ is ahead first, then $D_n$ is {\em ahead first}.

If in the crossing sequence of $\edge{r}{i}$ we have that $i+1,i+2,\ldots,n$ appear in this order, then $\edge{r}{i}$ is {\em ahead increasing}, and if they appear in descending order then it is {\em ahead decreasing}. If all the red edges are ahead increasing (respectively, decreasing), then the drawing $D_n$ itself is {\em ahead increasing} (respectively, {\em ahead decreasing}). If $D_n$ is ahead increasing or ahead decreasing, then it is {\em ahead monotone}.

Similarly, if in the crossing sequence of $\edge{r}{i}$ we have that $1,2,\ldots,i-1$ appear in this order, then $\edge{r}{i}$ is {\em behind increasing}, and if they appear in descending order then it is {\em behind decreasing}. If all the red edges are behind increasing (respectively, decreasing), then the drawing $D_n$ itself is {\em behind increasing} (respectively, decreasing). If $D_n$ is behind increasing or behind decreasing, then it is {\em behind monotone}.

Thus, for instance, if $D_n$ is ahead first, ahead decreasing, and behind increasing, then for each $i\in[n]$ the crossing sequence of the red edge $\edge{r}{i}$ is $n, n{-}1, \ldots, i{+}2,i{+}1, 1,2,\ldots,i{-}1$ (with obvious adjustments for $i=1$ and $i=n$).

If all the crossings of $r_i$ with $b_{i+1},b_{i+2},\ldots,b_n$ are positive (respectively, negative) then $r_i$ is {\em ahead positive} (respectively, {\em ahead negative}). If all the red edges are ahead positive (respectively, if they are all ahead negative) then $D_n$ itself is {\em ahead positive} (respectively, {\em ahead negative}). If $D_n$ is ahead positive or ahead negative, then it is {\em ahead consistent}.

Similarly, if all the crossings of $r_i$ with $b_{1},b_{2},\ldots,b_{i-1}$ are positive (respectively, negative) then $r_i$ is {\em behind positive} (respectively, {\em behind negative}). If all the red edges are behind positive (respectively, if they are all behind negative) then $D_n$ itself is {\em behind positive} (respectively, {\em behind negative}). If $D_n$ is behind positive or behind negative, then it is {\em behind consistent}.

Regarding blue edges, we need only pay attention to $b_1$ and $b_n$. If in $D_n$ the crossing sequence of $b_1$ is $2,\ldots,n-1,n$ (respectively, $n,n-1,\ldots,2$) then we say that $D_n$ is $b_1$-{\em increasing} (respectively, $b_1$-{\em decreasing}). In either case we say that $D_n$ is $b_1$-{\em monotone}. Similarly, if in $D_n$ the crossing sequence of $b_n$ is $1,2,\ldots,n-1$ (respectively, $n-1,\ldots,2,1$) then $D_n$ is $b_n$-{\em increasing} (respectively, $b_n$-{\em decreasing}). In either case we say that $D_n$ is $b_n$-{\em monotone}.

Let $\Sigma$ be a surface. An application of Ramsey's theorem for $3$-uniform hypergraphs shows that if there exist arbitrarily large $2$-drawings in $\Sigma$, then there exist arbitrarily $2$-drawings in $\Sigma$ that are ahead first. With successive applications of Ramsey's theorem we get a series of additional properties, ending up with the following.

\begin{proposition}\label{pro:2cr1}
Let $\Sigma$ be a surface. If there exist arbitrarily large $2$-drawings in $\Sigma$, then there exist arbitrarily large $2$-drawings in $\Sigma$ that are ahead first, ahead monotone, behind monotone, ahead consistent, behind consistent, $b_1$-monotone, and $b_n$-monotone.
\end{proposition}


\section{Proof of Lemma~\ref{lem:2cr}}\label{sec:2cr}

Lemma~\ref{lem:2cr} claims that no fixed surface $\Sigma$ can host $2$-drawings of $K_{2,m}$ for arbitrarily large values of $m$. In view of Proposition~\ref{pro:2cr1}, in order to prove the lemma it suffices to show that {\em no fixed surface can host $2$-drawings of $K_{2,n}$ that are ahead first, ahead monotone, behind monotone, ahead consistent, behind consistent, $b_1$-monotone, and $b_n$-monotone.}

Even though evidently the focus has been reduced to a finite number of possibilities, at first glance we seem to have a grueling task ahead. Indeed, since there are $2$ ways in which a $2$-drawing can be ahead monotone, $2$ ways in which it can be behind monotone, $2$ ways in which it can be ahead consistent, $2$ ways in which it can be behind consistent, $2$ ways in which it can be $b_1$-monotone, and $2$ ways in which it can be $b_n$-monotone, in principle we need to investigate $2^6=64$ different possibilities. Fortunately, as we shall see shortly, all cases are disposed of with a handful of simple arguments. 

Each of the $64$ cases is dealt with in one of the next three claims. For instance, the next statement swiftly takes care of $60$ cases with the same easy observation used to prove Claim~\ref{cla:1cla1}.
\newpage
\begin{claim}\label{cla:2cla1}
Let $D_n$ be a $2$-drawing of $K_{2,n}$ in a surface $\Sigma$, where $n$ is even. Suppose that $D_n$ is either:
\begin{description}[left=0.7cm]
\item[(i)] behind positive and $b_1$-increasing; or 
\item[(ii)] behind negative and $b_1$-decreasing; or 
\item[(iii)] ahead positive and $b_n$-increasing; or 
\item[(iv)] ahead negative and $b_n$-decreasing; or 
\item[(v)] ahead increasing, and ahead positive; or 
\item[(vi)] ahead decreasing, and ahead negative; or 
\item[(vii)] behind increasing, and behind positive; or 
\item[(viii)] behind decreasing, and behind negative.
\end{description}
Then $\genus(\Sigma)\ge (n-2)/2$.
\end{claim}

With the following claim, we will deal with $3$ of the $4$ remaining cases. 

\begin{claim}\label{cla:2cla2}
Let $D_n$ be a $2$-drawing of $K_{2,n}$ in a surface $\Sigma$. Suppose that $D_n$ is ahead first, and in addition it is either:
\begin{description}[left=0.7cm]
\item[(i)] ahead increasing, ahead negative, and $\edge{b}{n}$-increasing; or
\item[(ii)] behind increasing, behind negative, and $\edge{b}{1}$-increasing.
\end{description}
Then there exists a $1$-drawing of $K_{2,n-1}$ in $\Sigma$.
\end{claim}

Finally, with the next statement, we take care of the single remaining case.

\begin{claim}\label{cla:2cla3}
Let $D_n$ be a $2$-drawing of $K_{2,n}$ in a surface $\Sigma$, where $n\equiv 1 \inmod{3}$. Suppose that $D_n$ is of ahead first, ahead decreasing, behind decreasing, ahead positive, behind positive, and $\edge{b}{1}$-{decreasing}. Then $\genus(\Sigma)\ge (n-1)/3$.
\end{claim}

We defer the proofs of the claims for the moment, and show that they imply Lemma~\ref{lem:2cr}. 

\begin{proof}[Proof of Lemma~\ref{lem:2cr}, assuming Claims~\ref{cla:2cla1},~\ref{cla:2cla2}, and~\ref{cla:2cla3}]
Let $\Sigma$ be a surface. To prove Lemma~\ref{lem:2cr} we exhibit an integer $n$ and show that there cannot exist a $2$-drawing of $K_{2,n}$ in $\Sigma$.

To define $n$, first let $N$ be the maximum integer such that there is a $1$-drawing of $K_{2,N}$ in $\Sigma$. The existence of $N$ is guaranteed from Lemma~\ref{lem:1cr}. We let $n$ be any integer such that $n\ge\max\{3\,\genus(\Sigma)+3,N + 2\}$ and $n\equiv 1\inmod{3}$. 

By Proposition~\ref{pro:2cr1} there is a $2$-drawing $D_n$ of $K_{2,n}$ that is ahead first, ahead monotone, behind monotone, ahead consistent, behind consistent, $b_1$-monotone, and $b_n$-monotone.

A trivial analysis shows that $D_n$ either (a) satisfies one of (i)--(viii) in Claim~\ref{cla:2cla1}; or (b) it satisfies one of (i)--(ii) in Claim~\ref{cla:2cla2}; or (c) it satisfies the hypotheses in Claim~\ref{cla:2cla3}.

Suppose that (a) holds. Note that $n$ is even. By Claim~\ref{cla:2cla1}, we have that $\genus(\Sigma) \ge (n-2)/2$, that is, $n \le 2\,\genus(\Sigma)+2$. But this contradicts the choice of $n$, which in particular implies that $n \ge 2\,\genus(\Sigma) + 3$. 

If (b) holds, Claim~\ref{cla:2cla2} guarantees the existence of a $1$-drawing of $K_{2,n-1}$ in $\Sigma$, and from this it follows that $n-1 \le N$. But this contradicts the choice of $n$, which in particular implies that $n-1 \ge N+1$. 

Finally, suppose that (c) holds. By Claim~\ref{cla:2cla3}, we have that $\genus(\Sigma) \ge (n-1)/3$, that is, $n \le 3\,\genus(\Sigma)+1$. But this contradicts the choice of $n$, which in particular implies that $n \ge 3\,\genus(\Sigma)+2$.
\end{proof}

We emphasize that in the proofs of the claims we adopt once again the perspective used in the proof of Claims~\ref{cla:1cla1} and~\ref{cla:1cla2}, where crossings are regarded as degree $4$ vertices. Thus once again $\edge{b}{i}$ and $\edge{r}{i}$ are seen as paths, for $i=1,\ldots,n$.

\def\ric#1#2{{\bigl|_{\edge{#1}{#2}}}}

\begin{proof}[Proof of Claim~\ref{cla:2cla1}]
This claim follows using an argument totally analogous to the one used in the proof of Claim~\ref{cla:1cla1}, to find a restriction $J_n$ of $D_n$ with $2n+1$ vertices, $3n-2$ edges, and one facial walk.

For $i\in\{1,n\}$, we use $D_n \bigl|_{\edge{b}{i}}$ to denote the restriction of $D_n$ to the blue path $\edge{b}{i}$ and all the red paths with the exception of $\edge{r}{i}$. Similarly, $D_n \bigl|_{\edge{r}{i}}$ is the restriction of $D_n$ to the red path $\edge{r}{i}$ and all the blue paths with the exception of $\edge{b}{i}$. For instance, Figure~\ref{fig:4700} in Section~\ref{sec:1cr} illustrates $D_n\bigl|_{\edge{b}{1}}$ if $D_n$ is behind positive and $b_1$-increasing, that is, if $D_n$ is as in (i) in Claim~\ref{cla:2cla1}. 

As an additional example, in Figure~\ref{fig:4400} we illustrate $D_n\bigl|_{\edge{r}{n}}$ if $D_n$ is behind decreasing and behind negative, as in (viii) in Claim~\ref{cla:2cla1}.

\def\tf#1{{\Scale[2.5]{#1}}}
\def\tz#1{{\Scale[2.5]{#1}}}
\begin{figure}[ht!]
\centering
\hglue 0cm \scalebox{0.38}{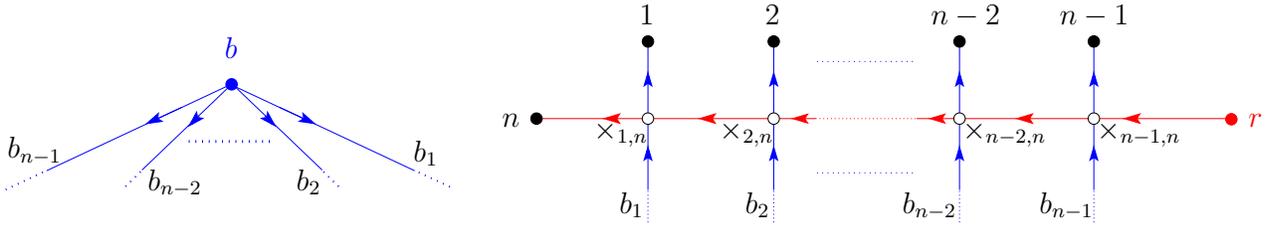}
\caption{Illustration of the proof of Claim~\ref{cla:2cla1}.}
\label{fig:4400}
\end{figure}
\def\tf#1{{\Scale[2.4]{#1}}}
\def\tz#1{{\Scale[2.0]{#1}}}

If (i) or (ii) holds, let $J_n:=D_n\ric{b}{1}$. If (iii) or (iv) holds, let $J_n:=D_n\ric{b}{n}$.  If (v) or (vi) holds, let $J_n:=D_n\ric{r}{1}$. Finally, if (vii) or (viii) holds, let $J_n:=D_n\ric{r}{n}$. In all these cases, as in the proof of Claim~\ref{cla:1cla1}, it is easily verified that $J_n$ has $2n+1$ vertices, $3n-2$ edges, and one facial walk. Using~\eqref{eq:genuseq}, we conclude that the genus of the host surface of $J_n$ (which is the host surface of $D_n$) is at least $(n-2)/2$. 
\end{proof}

\begin{proof}[Proof of Claim~\ref{cla:2cla2}]
Suppose that (i) holds. The idea of the proof is quite simple. First, we remove certain parts of $D_n$, thus obtaining a drawing $D_n'$. We then contract a path in $D_n'$ to a vertex, and show that as a result we obtain a $1$-drawing of $K_{2,n-1}$.

Recall that in a $2$-drawing of $K_{2,n}$, for each $i=1,\ldots,n$ the red path $\edge{r}{i}$ crosses the blue path $\edge{b}{j}$ for every $j\neq i$. In particular, for $i=1,2,\ldots,n-1$, the red path $\edge{r}{i}$ crosses the blue path $\edge{b}{n}$. As illustrated in Figure~\ref{fig:4800}(a), we let $s_i$ be the part of the red path $\edge{r}{i}$ from $r$ to its crossing $\cross{n}{i}$ with the blue path $\edge{b}{n}$, and we let $t_i$ be the rest of $\edge{r}{i}$, that is, the part from $\cross{n}{i}$ to vertex~$i$.  In this figure, the parts $t_i$ are drawn thick. 

\def\te#1{{\Scale[2.6]{#1}}}
\def\tf#1{{\Scale[2.4]{#1}}}
\def\tz#1{{\Scale[2.2]{#1}}}
\begin{figure}[ht!]
\centering
\hglue 0.1cm \scalebox{0.35}{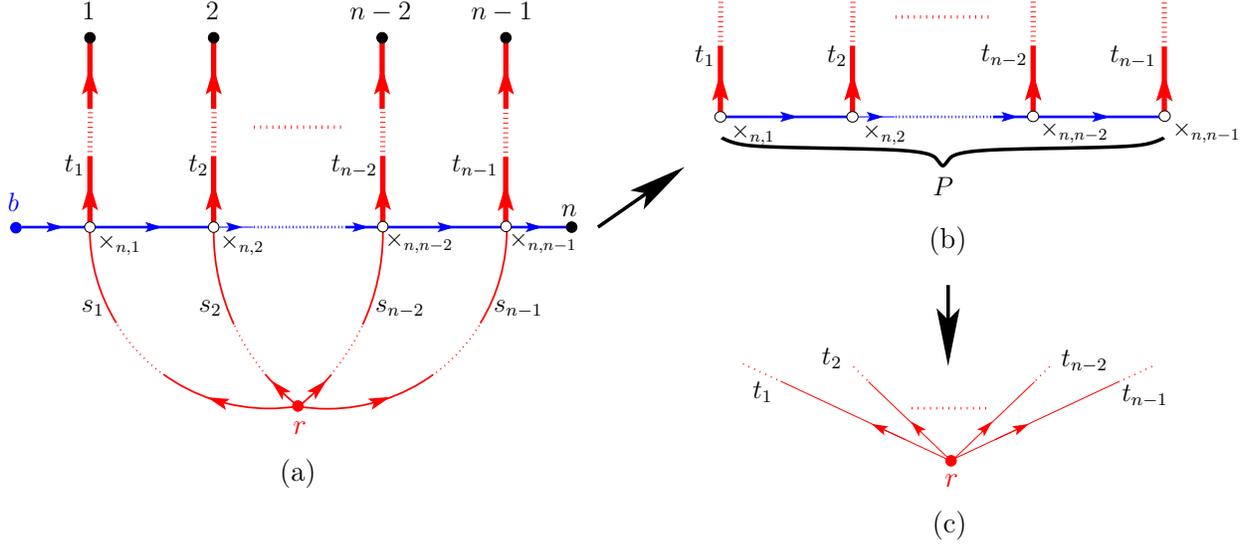}
\caption{Illustration of the proof of Claim~\ref{cla:2cla2}.}
\label{fig:4800}
\end{figure}
\def\te#1{{\Scale[2.2]{#1}}}
\def\tf#1{{\Scale[2.4]{#1}}}
\def\tz#1{{\Scale[2.0]{#1}}}

As we also illustrate in Figure~\ref{fig:4800}(a), the assumption that $D_n$ is $b_n$-increasing implies that as we traverse $\edge{b}{n}$ from $b$ to $n$ we find $\cross{n}{1},\cross{n}{2},\ldots,\cross{n}{n-1}$ in this order. Also note that the illustration in Figure~\ref{fig:4800}(a) reflects that $D_n$ is ahead negative, and so in particular $\cross{n}{1},\cross{n}{2},\ldots,\cross{n}{n-1}$ are negative crossings.

We now discard from $D_n$ the red path $\edge{r}{n}$, the red vertex $r$, and for $i{=}1,2,\ldots,n{-}1$ we discard the part $s_i$ of the red path $\edge{r}{i}$. Finally, we also discard the blue edge from ${b}$ to ${\cross{n}{1}}$ and also the blue edge from $\cross{n}{n-1}$ to $n$. We let $D_n'$ be the drawing obtained after removing these vertices and edges from $D_n$. In Figure~\ref{fig:4800}(b) we illustrate the part of $D_n'$ that consists of what remains of the blue path $\edge{b}{n}$ (which is the blue path $P$ from $\cross{n}{1}$ to $\cross{n}{n-1}$) and its incident red edges.

Finally, as illustrated in Figure~\ref{fig:4800}(c), we contract in $D_n'$ the whole blue path $P$ to a vertex that we label $r$. Note that there is no ambiguity in the use of this label, since the original red vertex $r$ was discarded at the beginning of the procedure. Let $F_{n-1}$ denote the drawing thus obtained. We conclude the proof of Case (i) by showing that $F_{n-1}$ is a $1$-drawing of $K_{2,n-1}$, where the blue paths are $\edge{b}{1},\edge{b}{2},\ldots,\edge{b}{n-1}$, and the red paths are $\edge{t}{1},\edge{t}{2},\ldots,\edge{t}{n-1}$.

To prove that $F_{n-1}$ is a $1$-drawing of $K_{2,n-1}$, we need to show that (I) the paths $t_1,t_2,\ldots,t_{n-1}$ leave $r$ in this clockwise cyclic order; and (II) for each pair of distinct $i,j\in\{1,2,\ldots,n-1\}$, there is exactly one crossing in the canonical cycle that consists of $\edge{b}{i},\edge{b}{j},\edge{t}{i}$, and $\edge{t}{j}$.

As illustrated in Figure~\ref{fig:4800}(c), (I) follows simply because the original drawing $D_n$ is $b_n$-increasing. Finally, to see that (II) holds it suffices to note that for $i=1,2,\ldots,n-1$, the red path $t_i$ crosses the blue paths $\edge{b}{1},\edge{b}{2},\ldots,\edge{b}{i-1}$, and it does not cross any of the blue paths $\edge{b}{i+1},\ldots,\edge{b}{n}$: this follows since the original drawing $D_n$ is ahead first and ahead increasing. Thus (I) and (II) hold, and so $F_{n-1}$ is a $1$-drawing of $K_{2,n-1}$ in $\Sigma$, as claimed.

If (ii) holds, we proceed in a totally analogous manner. In this case we discard the red vertex~$r$, the red path $\edge{r}{1}$, and for $i=2,\ldots,n$ the part of $\edge{r}{i}$ from $r$ to its crossing $\cross{1}{i}$ with the blue path $\edge{b}{1}$. Next, we discard the portions of $\edge{b}{1}$ from $b$ to $\cross{1}{2}$ and from $\cross{1}{n}$ to $1$. Finally, we contract the part of $\edge{b}{1}$ from $\cross{1}{2}$ to $\cross{1}{n}$ to a vertex we label $r$, and as a result we obtain a $1$-drawing of $K_{2,n-1}$ in $\Sigma$.
\end{proof}

\begin{proof}[Proof of Claim~\ref{cla:2cla3}]
As we illustrate in Figure~\ref{fig:7400}, the assumption that $D_n$ is ahead decreasing, ahead positive, behind decreasing, and behind positive, implies that for $i=2,\ldots,n-1$, as we traverse the red path $\edge{r}{i}$ from $r$ to $i$, we cross the blue paths $\edge{b}{n},\edge{b}{n-1},\ldots,\edge{b}{i+1},\edge{b}{i-1},\ldots,\edge{b}{2},\edge{b}{1}$ in this order, and we cross all these paths from their left-hand sides.

\def\te#1{{\Scale[2.2]{#1}}}
\def\tf#1{{\Scale[2.2]{#1}}}
\def\tz#1{{\Scale[2.2]{#1}}}
\begin{figure}[ht!]
\centering
\hglue 0 cm \scalebox{0.42}{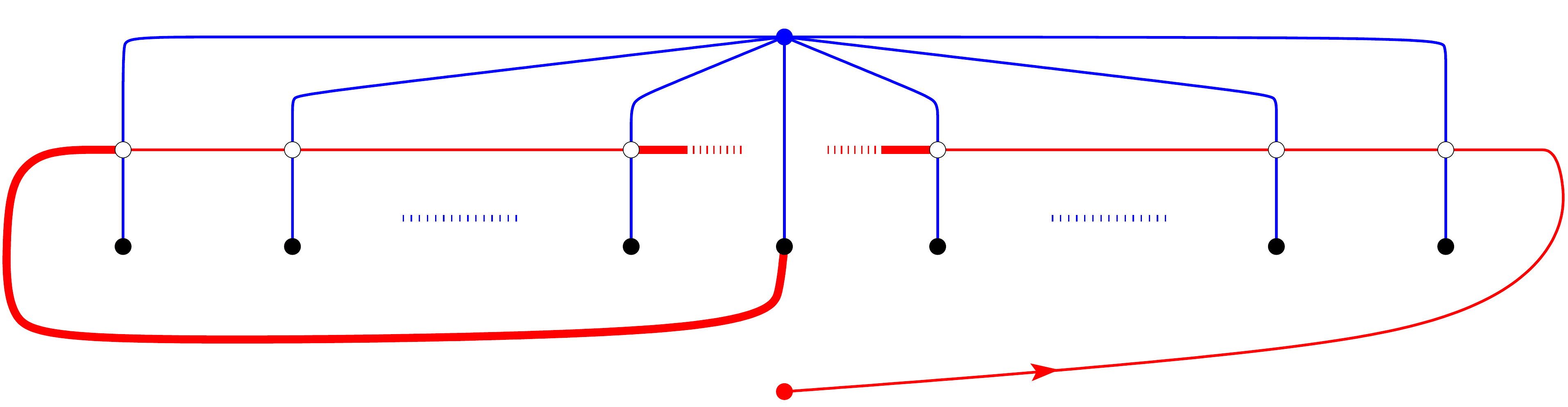}
\caption{Illustration of the proof of Claim~\ref{cla:2cla3}.}
\label{fig:7400}
\end{figure}
\def\tf#1{{\Scale[2.4]{#1}}}
\def\tz#1{{\Scale[2.0]{#1}}}

As we illustrate in this figure, after crossing $\edge{b}{i+1}$ at $\cross{i+1}{i}$, the path $\edge{r}{i}$ reaches $\edge{b}{i-1}$ at $\cross{i-1}{i}$. Thus there is a red edge with endpoints $\cross{i+1}{i}$ and $\cross{i-1}{i}$; we use $s_i$ to denote this red edge. We also note that as we traverse $\edge{r}{i}$, the last crossing we encounter is its crossing $\cross{1}{i}$ with $\edge{b}{1}$, and after this crossing we arrive in vertex $i$. Thus there is a red edge with endpoints $\cross{1}{i}$ and $i$; we use $t_i$ to denote this red edge. Both $s_i$ and $t_i$ are thick in Figure~\ref{fig:7400}.

We let $J_n$ be the restriction of $D_n$ to the part that consists of the blue paths $\edge{b}{1},\ldots,\edge{b}{n}$ plus the segments $s_i$ and $t_i$ for $i=3,6,9,\ldots,n-1$. See Figure~\ref{fig:7700}. {As we illustrate in that figure, since $D_n$ is $b_1$-decreasing it follows that the crossing-vertices $\cross{1}{n-1},\ldots,\cross{1}{6},\cross{1}{3}$ appear in this order as we traverse $\edge{b}{1}$ from $b$ to $1$.}

\def\te#1{{\Scale[2.4]{#1}}}
\def\tf#1{{\Scale[2.4]{#1}}}
\def\tz#1{{\Scale[2.4]{#1}}}
\begin{figure}[ht!]
\centering
\hglue -1.4 cm \scalebox{0.38}{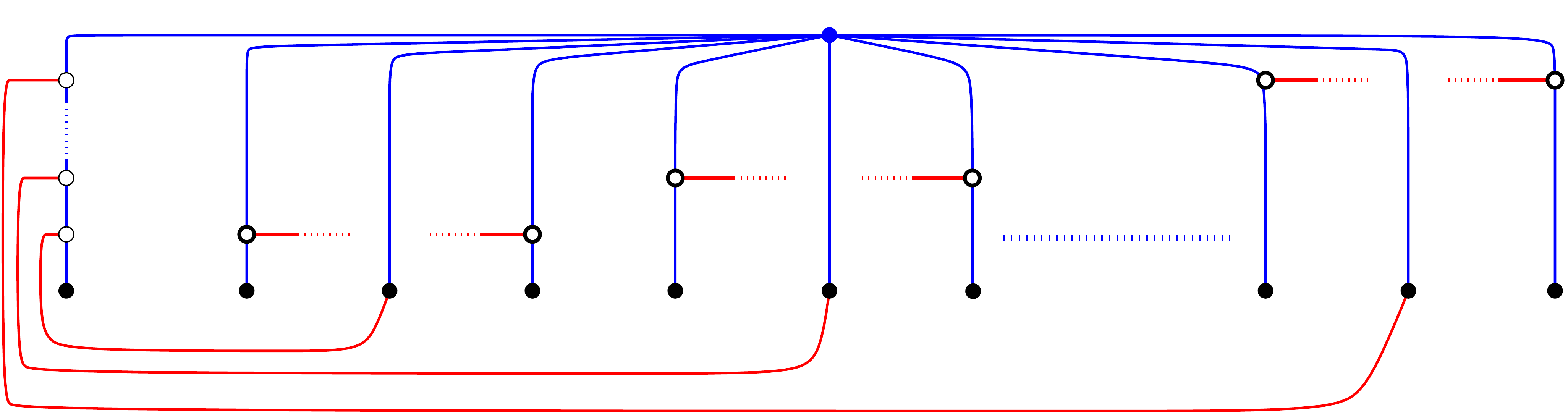}
\caption{Illustration of the proof of Claim~\ref{cla:2cla3}.}
\label{fig:7700}
\end{figure}
\def\tf#1{{\Scale[2.4]{#1}}}
\def\tz#1{{\Scale[2.0]{#1}}}

We note that $J_n$ has $2n$ vertices and $8n/3 -5/3$ edges, and it is straightforward to verify that it has exactly one facial walk. Using~\eqref{eq:genuseq} we obtain that the genus of the host surface of $J_n$ (which is the host surface of $D_n$) is at least $(n-1)/3$.
\end{proof}

\section{Concluding remarks and open questions}\label{sec:con}

It is well known (see for instance~\cite{grosstucker,moharthomassen}) that every rotation system is realized by an embedding on a surface. Since an embedding is an (particularly nice) instance of a simple drawing, it immediately follows that every rotation system, and in particular every complete rotation system, is simply realizable in every surface whose genus is sufficiently large.

Archdeacon defined the {\em crossing genus} $\croge(\Pi)$ of a rotation system $\Pi$ as the smallest genus of a surface in which $\Pi$ can be simply realized. He put forward the following.

\vglue 0.2 cm
\noindent{\bf Problem} (Archdeacon~\cite{Arch1}). {\em Given a complete rotation system, calculate its crossing genus.}
\vglue 0.2 cm

For each nonnegative integer $g$, let $f(g)$ be the least integer such that there exists a complete rotation system of size $f(g)$ that is not simply realizable in the surface of genus $g$. For instance, $f(0)=4$, since there exist complete rotation systems of size $4$ that are not simply realizable in the sphere (equivalently, the plane), but trivially every complete rotation system of size $3$ is simply realizable in the sphere.

Theorem~\ref{thm:main} implies that $f(g)$ is well-defined. On the other hand, in our arguments we use Ramsey's theorem repeatedly, and so the upper bound we can prove for $f(g)$ is multiply exponential in $g$. Can one prove a remarkably better upper bound for $f(g)$? For instance, is it true that $f(g) = O(2^g)$? 

How about lower bounds for $f(g)$? Any rotation system of size $n$ can be realized as a 2-cell embedding (and hence as a simple drawing) of $K_n$ in some surface. If the rotation system defines $\phi$ facial walks, then the genus of this surface is $(1/2)(2-|V(K_n)|+|E(K_n)|-\phi)\le(1/2)(1-n+\binom{n}{2})=(1/4)(n-1)(n-2) < n^2/4$. Therefore every rotation system of size $n$ can be simply realized in the surface of genus $n^2/4$, and this implies that $f(g) =\Omega(\sqrt{g})$. Can one prove a significantly better lower bound for $f(g)$? For instance, is it true that $f(g) = \Omega(g)$?

\section*{Acknowledgements}
\begin{minipage}[t]{0.75\textwidth}
	This work was supported by the H2020-MSCA-RISE project 734922-CONNECT. Rosna Paul and Alexandra Weinberger acknowledge the support of the Austrian Science Fund (FWF): W1230. Gelasio Salazar acknowledges the support of CONACYT under Proyecto Ciencia de Frontera 191952.
	\end{minipage}\hfill
\begin{minipage}[t]{0.25\textwidth}
	\vspace*{-2.35cm}
	\hspace*{-1.7cm}
	{\includegraphics[scale=0.28]{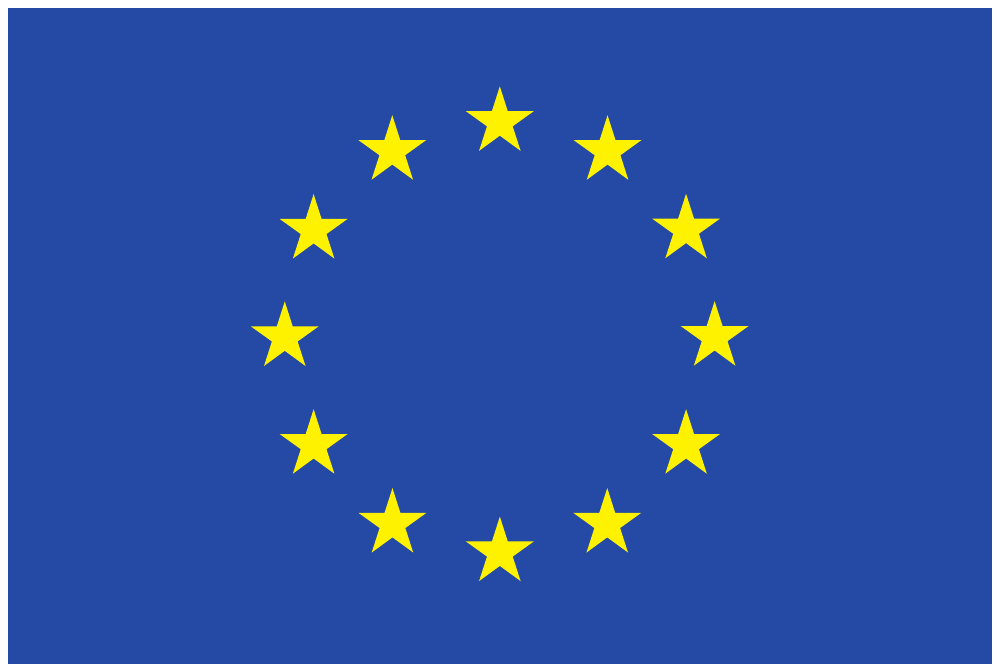}}
\end{minipage}
\bibliographystyle{abbrv}
\bibliography{rotsys.bib}

\begin{thebibliography}{1}

\bibitem{EuroCG15}
B.~\'Abrego, O.~Aichholzer, S.~Fern\'andez-Merchant, T.~Hackl, J.~Pammer,
  A.~Pilz, P.~Ramos, G.~Salazar, and B.~Vogtenhuber.
\newblock {{All Good Drawings of Small Complete Graphs}}.
\newblock In {\em Proc. $31^{st}$ European Workshop on Computational Geometry
  EuroCG '15}, pages 57--60, Ljubljana, Slovenia, 2015.

\bibitem{Arch1}
D.~Archdeacon.
\newblock Problems in topological graph theory ---questions {I} can't
  answer---.
\newblock {\em Yokohama Mathematical Journal}, 47:89--92, 1999.

\bibitem{cardinalfelsner}
J.~Cardinal and S.~Felsner.
\newblock Topological drawings of complete bipartite graphs.
\newblock {\em J. Comput. Geom.}, 9(1):213--246, 2018.

\bibitem{grosstucker}
J.~L. Gross and T.~W. Tucker.
\newblock {\em Topological graph theory}.
\newblock Wiley-Interscience Series in Discrete Mathematics and Optimization.
  John Wiley \& Sons Inc., New York, 1987.

\bibitem{Kyncl1}
J.~Kyn\v{c}l.
\newblock Simple realizability of complete abstract topological graphs in {{\bf
  {P}}}.
\newblock {\em Discrete and Computational Geometry}, 45:383--399, 2011.

\bibitem{Kyncl2}
J.~Kyn\v{c}l.
\newblock Simple realizability of complete abstract topological graphs
  simplified.
\newblock {\em Discrete and Computational Geometry}, 64:1--27, 2020.

\bibitem{moharthomassen}
B.~Mohar and C.~Thomassen.
\newblock {\em Graphs on Surfaces}.
\newblock Johns Hopkins series in the mathematical sciences. Johns Hopkins
  University Press, 2001.

\end{thebibliography}

\end{document}